\documentclass[12pt]{article}
\usepackage[english,russian]{babel}
\begin{document}
\centerline{\Large Remark on the Smale's Problem 9}

\bigskip
\centerline{\Large Neil K. Bakirov }

\bigskip
\centerline{Institution of Russian Academy of Sciences }

\centerline{Institute of Mathematics with Computing Centre}

\bigskip
\begin{abstract}
We consider the S.Smale's 9th problem on feasibility of the linear
system of inequalities in connection with a linear programming
problem.
\end{abstract}


\bigskip\bigskip
In [1] it was posed the following problem 9. Let $A$ be an $m\times
n$ real matrix and $b\in R^m$.

\bigskip
\noindent{\it Problem 9: Is there a polynomial time algorithm over
the real numbers which decides the feasibility of the linear system
of inequalities $Ax\ge b$?}

\bigskip
The answer is -- yes there is.

It may be shown by reduction of the problem 9 to the linear
programming problem with certainly feasible conditions. So the
problem 9 has a required polynomial time algorithm as the linear
programming problem has such an algorithm (L.Khachiyan, 1979 and
subsequently N.Karmarkar, 1984).

{\bf Reduction:} Let
$$
\Delta b\in R^m,\qquad c=(1,1,...,1)\in R^m,\qquad F(x,\Delta
b)=(c,\Delta b),
$$
where $(.,.)$ -- is the inner product, $c$ -- is the vector with
unit coordinates.

Consider  the following linear programming problem:
$$
\min_{x\in R^n}F(x,\Delta b)\hspace{1cm}\mbox{subject
to}\hspace{1cm} Ax\ge b-\Delta b,\quad \Delta b\ge0\eqno(1)
$$
here the goal function is bounded from below $F(x,\Delta b)\ge0$ and
conditions $Ax\ge b-\Delta b,\Delta b\ge0$ are certainly feasible.

\bigskip\bigskip
{\bf Theorem.} The linear system of inequalities $Ax\ge b$ is
feasible if and only if in the problem (1)
$$
\min_{x\in R^n}F(x,\Delta b)=0.\eqno(2)
$$
Proof:

Necessity. If there exist such a vector $x$ for which $Ax\ge b$ then
(2) is fulfilled with $\Delta b=0.$

Sufficiency. If (2) is fulfilled for some $(x,\Delta b)$ then
$\Delta b=0$ and therefore $x$ is a solution of the linear system of
inequalities $Ax\ge b$.

\bigskip
Theorem is proved.

\bigskip
{\bf Corollary.} To check solvability of the linear system of
inequalities $Ax\ge b$ one can first solve the problem (1) using a
polynomial time algorithm and then check the condition (2) which
solves the problem 9.

\bigskip
\centerline{\bf Literature}

\bigskip
1. Smale S., Mathematical Problems for the Next Century,
Mathematics: frontiers and perspectives, eds. Arnold, V., Atiyah,
M., Lax, P. and Mazur, B., Amer. Math. Soc., 2000.

2. http://en.wikipedia.org/wiki/Linear\_programming

\bigskip\bigskip
\noindent Neil K. Bakirov,

Institute of Mathematics of Ufa Scientific Centre

of the Russian Academy of Sciences,

450000, Chernyshevsky str., Ufa, Russia,

bakirovnk@rambler.ru
\end{document}